\documentclass[12pt]{article}
\usepackage{amssymb}
\usepackage{amsmath}
\usepackage{latexsym}
\usepackage{mathrsfs}
\usepackage{graphicx}

\numberwithin{equation}{section}
\newtheorem{thm}[equation]{Theorem}

\newtheorem{defn}[equation]{Definition}
\newtheorem{rem}[equation]{Remark}
\newtheorem{lem}[equation]{Lemma}
\newtheorem{corol}[equation]{Corollary}

\title{On the sharpness of a certain spectral stability estimate for the Dirichlet Laplacian\footnote{Published in {\it Eurasian Mathematical Journal, } Volume 1, Number 1 (2010), 111-122. Available also at the web page http://www.enu.kz/en/emj.php}}
\author{Pier Domenico Lamberti and Marco Perin }

\date{ \today }

\begin{document}

\newcommand{\rea}{\mathbb{R}}

\maketitle

%
%
%

\noindent
{\bf Abstract:} We consider a spectral stability estimate by Burenkov and Lamberti concerning the variation of the eigenvalues of second order uniformly elliptic operators on variable open sets in the $N$-dimensional euclidean space, and  we prove that it is sharp for any dimension $N$. This is done by studying  the eigenvalue problem for the Dirichlet Laplacian on special open sets inscribed in suitable spherical cones.

\vspace{11pt}

\noindent
{\bf Keywords:} Elliptic equations, Dirichlet boundary conditions,
stability of eigenvalues, sharp estimates,
domain perturbation.

\vspace{6pt}
\noindent
{\bf 2000 Mathematics Subject Classification:} 35P15, 35J40, 47A75, 47B25.
%
%



\section{Introduction}

Let $\Omega $ be a bounded open set in ${\mathbb{R}}^N$. We consider the eigenvalue problem for the Dirichlet Laplacian
\begin{equation}\label{mainpr}
\left\{\begin{array}{ll}-\Delta u=\lambda u, & \ {\rm in }\ \Omega ,\\
u=0, & \ {\rm on }\ \partial \Omega .
\end{array}  \right.
\end{equation}
As is well known, problem (\ref{mainpr}) has a non-decreasing  divergent sequence of positive eigenvalues of finite multiplicity
$$
0<\lambda_1[\Omega]\le \lambda_2[\Omega]\le \dots \le \lambda_n[\Omega ] \le \dots
$$
Here each eigenvalue is repeated according to its  multiplicity. The dependence of $\lambda_n[\Omega ]$ on $\Omega $ has been studied by many authors. We refer to Burenkov, Lamberti and Lanza de Cristoforis~\cite{burlamlan} for a survey paper on this topic. In particular, several  papers have been devoted to the problem of finding explicit estimates
for the variation of $\lambda _n[\Omega ]$ upon variation of $\Omega $, see {\it e.g.,} Davies~\cite{da2000}, Burenkov and Lamberti~\cite{buladir, bulahigh}.

As a consequence of a general result proved in Burenkov and Lamberti~\cite[Cor.~5.16]{buladir}  for second order uniformly elliptic operators, we have that if $\Omega_1$ is a fixed open set and $p\in ]2, \infty ]$ is such that all the eigenfunctions $u $ of the Dirichlet Laplacian
on $\Omega _1$ satisfy the condition
\begin{equation}\label{sumcond}
\nabla u \in L^p(\Omega_1),
\end{equation}
then for each $n\in {\mathbb{N}}$ there exists $c_n>0$ such that
\begin{equation}
\label{bulaest}
|\lambda _n[\Omega _1]- \lambda _n[\Omega _2]|\le c_n |\Omega _1\setminus \Omega _2|^{1-\frac{2}{p}},
\end{equation}
for all open sets $\Omega_2\subset \Omega_1$ satisfying $|\Omega_1\setminus \Omega_2|<c_n^{-1}$. This result is valid  for open sets  $\Omega_1$ and $\Omega_2$ belonging to a uniform class  $C^{0,1}_M({\mathcal{A}})$ of open sets with Lipschitz continuous boundaries. Here  ${\mathcal{A}}$ denotes the `atlas' by the help of which the open sets are locally represented as subgraphs of Lipschitz continuous functions and $M>0$ is the uniform upper bound for the Lipschitz constants of  these functions. See Definition~\ref{class}, Theorem~\ref{bulathm} and Corollary~\ref{dimecor} for the precise statements.

It is clearly of interest to know whether the exponent in the right-hand side of (\ref{bulaest}) is sharp. An example confirming that  the exponent is sharp when  $N=2$ has been given in \cite{bulahigh}.
In this paper we develop an idea used in \cite{bulahigh} and we prove the sharpness of the exponent $1-2/p$ for any dimension $N$.

To do so, we consider a spherical cone $\Omega_{\beta }$ of angle $\beta \in [0, \pi [$ and a perturbation $\Omega_{\beta }(\epsilon)$, $\epsilon >0$ obtained by removing a ball of radius $\epsilon $ centered at the vertex of the cone.
In this case, we find the exact asymptotic behavior of the eigenvalues $\lambda [\Omega _{\beta }(\epsilon )]$ as $\epsilon \to 0$ and the range of exponents $p$ for which condition (\ref{sumcond}) is satisfied in $\Omega_{\beta }$.
These explicit computations have been performed in Section 2 and  confirm the sharpness of estimate (\ref{bulaest}) as explained in Remark \ref{mainrem}.

We note that the open sets $\Omega_{\beta }(\epsilon )$ have Lipschitz continuous boundaries. However, if $\beta >\pi /2$  the open sets $ \Omega_{\beta }(\epsilon )$ are not of class $C^{0,1}_M({\mathcal{A}})$ for the same fixed atlas ${\mathcal{A}}$. Thus, in order to prove that the result of Theorem \ref{bulathm} is in fact sharp, we have to refine the construction by replacing the open sets $\Omega_{\beta }(\epsilon )$ by suitable open sets $\tilde \Omega _{\beta }(\epsilon )$ which enjoy the required properties. This construction is done in Section 3 where it is finally  proved that the result of Theorem \ref{bulathm} is sharp for any dimension $N$ (see Theorem \ref{shathm}).

\section{An example}

In this section we consider the eigenvalue problem (\ref{mainpr}) on a spherical cone $\Omega _{\beta }$ in ${\mathbb{R}}^N$ with angle $\beta $, and we study the variation of the eigenvalues upon suitable deformations $\Omega _{\beta }(\epsilon )$ of $\Omega _{\beta }$ depending on the scalar parameter $\epsilon >0$. As mentioned in the introduction, the open sets $\Omega_{\beta }(\epsilon )$ are obtained by removing from $\Omega_{\beta }$ a ball of radius $\epsilon $ centered at the vertex of the cone.

To do so, it is convenient to use the spherical coordinates $(r,\theta _1, \dots , \theta_{N-1})$ and the corresponding transformation of coordinates
\begin{align*}
  x_1  &= r\cos \theta _1  \hfill \\
  x_{2}  &= r\sin \theta _1 \cos \theta _2  \hfill \\
  &\vdots  \hfill \\
  x_{N-1}  &= r\sin \theta _1 \sin \theta _2  \cdots \sin \theta _{N - 2} \cos \theta _{N - 1}  \hfill \\
  x_N  &= r\sin \theta _1 \sin \theta _2  \cdots \sin \theta _{N - 2} \sin \theta _{N - 1}  ,
\end{align*}
 with $\theta _1 , \ldots \theta _{N - 2}  \in [ {0,\pi } ]$, $\theta _{N - 1}  \in [ {0,2\pi } [$ (here it is understood that $\theta _1 \in [0, 2\pi [ $ if $N=2$).

If $N=2$ we set
\begin{equation}
\label{omegaN2}
\Omega _{\beta }=\biggl\{(r, \theta_1):\ r\in]0,1[,\ \ \theta_1\in [0,\beta [\,\cup\, ] 2\pi -\beta , 2\pi [ \biggr\} ,
\end{equation}
for all $\beta \in ]0, \pi [ $.
If $N\geq 3$ we set
\begin{equation}
\label{omegaN34}
\Omega _{\beta }=\{(r, \theta_1, \dots , \theta_{N-1}):\ r\in]0,1[,\ \ \theta_1 \in [0, \beta [\ \},
\end{equation}
for all $\beta \in ]0, \pi [ $.

In both cases, $\beta $ is  fixed and the perturbation $\Omega _{\beta }(\epsilon )$ of $\Omega $ under consideration is defined by
\begin{equation}
\label{omegaeps}
\Omega _{\beta } (\epsilon )= \Omega_{\beta } \cap \{ (r, \theta_1, \dots , \theta_{N-1}):\ r\in ]\epsilon ,1[\  \},
\end{equation}
for all $\epsilon\in ]0,1[$. It is also convenient to set $\Omega_{\beta }(0)=\Omega_{\beta}$.

Now we recall some  well-known facts in order to fix the notation for the sequel, see {\it e.g.} Kozlov, Maz'ya and Rossmann~\cite[\S 2.2]{koz}. Recall that the Laplace operator can be written in spherical coordinates as
$$
 \Delta f = \frac{\partial^2 f}{\partial r^2} + \frac{N-1}{r} \frac{\partial f}{\partial r} + \frac{1}{r^2} \delta f,
$$
where
$$
\delta  = \sum\limits_{j = 1}^{N - 1} {\frac{1}
{{q_j \left( {\sin \theta _j } \right)^{N - j - 1} }}\frac{\partial }
{{\partial \theta _j }}\left( {\left( {\sin \theta _j } \right)^{N - j - 1} \frac{\partial }
{{\partial \theta _j }}} \right)}
$$
is the Laplace-Beltrami operator on the unit sphere $S^{N-1}$ of ${\mathbb{R}}^N$ and
\[
\begin{array}{*{20}c}
   {q_1  = 1,} & {q_j  = \left( {\sin \theta _1 \sin \theta _2  \cdots \sin \theta _{j - 1} } \right)^2 ,} & {j = 2, \ldots ,N - 1}  \\
 \end{array} .
\]
 As is well known, problem (\ref{mainpr}) can be solved by separation of variables. Namely, by setting $u(r , \theta_1, \dots , \theta _{N-1})= R(r)\Theta (\theta_1, \dots , \theta_{N-1})$ and using $l(l+N-2 )$ as separation constant, one reduces the study of equation $-\Delta u=\lambda u$ to the study of the equations

\begin{equation}
\label{RadN}
r^2 \frac{\partial^2 R}{\partial r^2} + r \left(N - 1 \right) \frac{\partial R}{\partial r}  + \lambda R r^2  - l\left(l+N-2\right) R = 0,
\end{equation}
and
\begin{equation}
\label{TrasvN}
-\delta \Theta = l\left(l+N-2\right) \Theta \,.
\end{equation}
The function $\Theta $ is defined on the spherical cap
\begin{equation}
C_{\beta }= \overline{\Omega }_{\beta }\cap S^{N-1}
\end{equation}
and has to satisfy the Dirichlet boundary condition
\begin{equation}
\label{bbel}
\Theta =0 , \ \ \ {\rm on}\ \ \partial C_{\beta },
\end{equation}
where $\partial C_{\beta }$ denotes the boundary of $C_{\beta }$ in the sphere $S^{N-1 }$  (here $\overline{\Omega }_{\beta }$ denotes the closure
of $\Omega_{\beta }$ in ${\mathbb{R}}^N$).

We denote by  $\Sigma_{\beta }$ the set of positive real numbers $l$ such that there exists a nontrivial solution to equation (\ref{TrasvN}) satisfying the boundary condition (\ref{bbel}). As is known $\Sigma _{\beta }$  consists of an increasing sequence of positive real numbers. We set
$$
l_{\beta }=\min \Sigma_{\beta}
$$
and we recall that $l_{\beta} $ is the so-called {\it characteristic value of the spherical cap $C_{\beta }$}, see {\it e.g.,} C\'{a}mera~\cite{camera}.

The following lemma is probably well known. For the convenience of the reader, we include a proof.

\begin{lem}\label{charcap} If $\beta \in ]0,\pi /2[$ then $l_{\beta }>1$. If $\beta =\pi /2$ then $l_{\beta }=1$. If $\beta \in ]\pi /2, \pi [$ then $l_{\beta }<1 $.
\end{lem}

{\bf Proof.}  By the monotonicity of the eigenvalues of the Laplace-Beltrami operator with Dirichlet boundary conditions it follows that
if $\beta_1<\beta_2$ then $l_{\beta_1}>l_{\beta_2}$. Thus, in order to prove the lemma it is enough to prove that $l_{\pi /2 }=1$.
To do so, we recall that  $l_{\beta } $ is the smallest positive zero of the function $l \mapsto P^{-\mu }_{l+\mu}(\cos \beta )$ with $\mu =(N-3)/2 $,
where $P^{-\mu }_{l+\mu }$ denotes the Legendre function of the first kind, order $-\mu $ and degree  $l+\mu $, see C\'{a}mera~\cite[p.~75]{camera}. If $\beta =\pi /2$ we have to study the function $l\mapsto  P^{-\mu }_{l+\mu}(0)$. It is well known that
$$
 P^{-\mu }_{l+\mu}(0)= 2^{-\mu} \pi^{-1/2}   \cos \left(\frac{l\pi }{2}\right)\Gamma \left( \frac{1+l}{2}  \right)/ \Gamma \left(1+\frac{l}{2}+ \mu \right),
$$
(see {\it e.g.,} Abramowitz and Stegun~\cite[p.~334]{olv}) hence the smallest positive zero of  $P^{-\mu }_{l+\mu}(0)$ is $l=1$.\hfill $\Box$\\

\begin{rem}\label{camerarem} If $N=2$ it is straightforward to prove that $l_{\beta }=\pi /(2 \beta )$. It is also known that $l_{\beta }=(\pi -\beta )/\beta $ if $N=4$.
Moreover, the asymptotic behavior of $l_{\beta }$ as $\beta \to \pi $ is known. Namely, as $\beta \to \pi $ we have
$$ \begin{array}{ll}
l_{\beta }\sim  \frac{1}{2\log \frac{2}{\pi -\beta  }}  & {\rm if }\ $N=3$,\vspace{0.3cm}\\
l_{\beta }\sim \frac{\Gamma (N-2)}{\Gamma\left(\frac{N-1}{2}\right) \Gamma \left(\frac{N-3}{2} \right)} \left( \frac{\pi -\beta }{2}  \right)^{N-3}& {\rm if}\ N\geq 4\, .
\end{array}$$
We refer to C\'{a}mera~\cite{camera} for details and proofs.
\end{rem}

For $N=2$ the following result can be found {\it e.g.,} in Davies~\cite[p.~132]{da}.

\begin{thm}\label{thmsum} If $\beta \in ]0 , \pi /2]$ then the gradients of the eigenfunctions of the Dirichlet Laplacian on  $\Omega_{\beta }$ belong $L^{\infty }(\Omega_{\beta} )$; if $\beta \in ] \pi /2 , \pi [$ they belong
to $L^{p}(\Omega_{\beta })$ for all $p\in [2, N/(1-l_{\beta })[$. Moreover, if $\beta \in ] \pi /2 , \pi [$ then the gradient of any nonzero eigenfunction corresponding to the first eigenvalue does not belong to $L^{p}(\Omega_{\beta })$ with $p=N/(1-l_{\beta })$.
\end{thm}
{\bf Proof. } The eigenfunctions of the Dirichlet Laplacian on $\Omega_{\beta } $ are of the form $R(r )\Theta (\theta_1, \dots ,\theta_{N-1} )$ where, for some $l>0$,  $R$ satisfies equation (\ref{RadN}) subject to the boundary condition $R(1)=0$ and $\Theta $ satisfies (\ref{TrasvN}) subject to condition (\ref{bbel})
By setting
$
R\left( r \right) = {r^{1-\frac{N}{2}}}u\left( r \right)
$
and
$$
\nu_{l}=\frac{(N + 2 l - 2)}{2} ,
$$
it follows that $u$ satisfies the Bessel equation
\begin{equation}
\frac{\partial^2 u}{\partial r^2}+\frac{1}{r}\frac{\partial u}{\partial r} +\left(\lambda-\frac{\nu_l^2}{r^2}  \right) u=0.
\end{equation}
Taking into account that the eigenfunctions of the Dirichlet Laplacian are well known to be bounded, it follows that $u$ is a multiple of the function
\begin{equation}
r^{1-\frac{N}{2}} J_{\nu_l}\left( r\sqrt{\lambda }\right),
\end{equation}
where $J_{\nu_l }$ denotes the Bessel function of the first kind and order $\nu_l$.
Moreover, it follows that
\begin{equation}
\label{SolBesRN}
\frac{\partial R}{\partial r}  \sim r^{\nu_l - \frac{N}
{2}  }  = r^{l - 1},\ \ \ {\rm as}\ r\to 0 .
\end{equation}

The proof easily follows by Lemma~\ref{charcap}, formula (\ref{SolBesRN}), by observing that $\Theta$ is a smooth function and that
the  eigenfunctions corresponding to the first eigenvalue are obtained for $l=l_{\beta }$  (this can be proved by recalling that the eigenfunctions of the Laplace operator corresponding to the first eigenvalue are the only eigenfunctions which do not change sign and by observing that if $l=l_{\beta }$ then any solution to problem (\ref{TrasvN}) subject to the boundary conditions (\ref{bbel})  does not change sign in $C_{\beta }$).
\hfill $\Box $\\

\begin{rem} In Theorem~\ref{thmsum} it is in fact proved that if $u=R\Theta $ is an eigenfunction of the Dirichlet Laplacian in $\Omega_{\beta }$ with $\Theta $ satisfying equation (\ref{TrasvN}) for some $l>0$, then $\nabla u\in L^{\infty }(\Omega_{\beta } )$ if $l\geq 1$ and $\nabla u \in L^{p}(\Omega_{\beta })$ for all $p\in [2, N/(1-l)[$ if $l <1$.
\end{rem}

Given $l\in \Sigma_{\beta}$, in order to solve  problem (\ref{mainpr}) on the open set $\Omega _{\beta }(\epsilon )$, one has to solve equation (\ref{RadN}) subject to  the boundary conditions
$R(\epsilon )=R(1)=0$. This leads to the linear  system
\[
\left\{ \begin{array}{ll}
  C_1 J_{\nu_l}  ( {\sqrt \lambda  } ) + C_2 Y_{\nu_l} ( {\sqrt \lambda  } )   &= 0 \vspace{0.1cm}\\
  {C_1 J_{\nu_l} ( {\epsilon \sqrt \lambda  } ) + C_2 Y_{\nu_l}  ( {\epsilon \sqrt \lambda  } )} & = 0
\end{array}  \right.
\]
in the unknowns $C_1,C_2\in {\mathbb{R}}$, which admits nontrivial  solutions  if and only if
\begin{equation}
\label{EqCaratt}
J_{\nu_l} ( {\sqrt \lambda  } )Y_{\nu_l} ( {\epsilon \sqrt \lambda  } ) - Y_{\nu_l} ( {\sqrt \lambda  } )J_{\nu_l} ( {\epsilon \sqrt \lambda  } ) = 0.
\end{equation}
As is well known, the set of all the eigenvalues of the Dirichlet Laplacian on $\Omega _{\beta }(\epsilon )$ is given by the union of the sets of zeros of the cross-product equations (\ref{EqCaratt}). On the other hand, by similar considerations the eigenvalues of the Dirichlet Laplacian on $\Omega _{\beta }$ are given by the zeros of the equations
\begin{equation}
\label{EqCarattunp}
J_{\nu_l} (\sqrt{\lambda })=0.
\end{equation}

\begin{thm}\label{mainthm}
Let $\beta \in ]0,\pi [$. Let $l\in \Sigma _{\beta }$
and $\lambda_*\ne 0$ be a solution of equation (\ref{EqCarattunp}). Then for any $\epsilon >0$ sufficiently small  there exists
a unique solution $\lambda (\epsilon )$ to equation (\ref{EqCaratt}) such that $\lambda (\epsilon )\to \lambda_*$ as $\epsilon \to 0$ and it satisfies
\begin{equation}\label{main}
\lambda (\epsilon )= \lambda _*+ a|\Omega_{\beta } \setminus \Omega_{\beta }(\epsilon )|^{\frac{N+2l-2}{N}}+o(|\Omega_{\beta } \setminus \Omega_{\beta }(\epsilon )|^{\frac{N+2l-2}{N}}  ),\ \ {\rm as}\ \epsilon \to 0^+,
\end{equation}
for some constant $a> 0$.
In particular
\begin{equation}
\label{main2}
\lambda _1 [\Omega _{\beta }(\epsilon )]= \lambda_1[\Omega _{\beta }]+ b|\Omega_{\beta } \setminus \Omega_{\beta }(\epsilon )|^{\frac{N+2l_{\beta } -2}{N}}+o\biggl(|\Omega_{\beta } \setminus \Omega_{\beta }(\epsilon )|^{\frac{N+2l_{\beta }-2}{N}}  \biggr),\ \ {\rm as}\ \epsilon \to 0^+,
\end{equation}
for some constant $b> 0$.
\end{thm}

{\bf Proof.} Note that $2\nu _l>1$ and
\begin{equation}
\label{ratio}
\frac{J_{\nu_{l}} (t)}{Y_{\nu_l}(t)}= t^{2\nu_l }h(t),
\end{equation}
for all $t>0$ sufficiently small, where $h$ is a function of class $C^1$ defined in a suitable neighborhood of zero in $[0,\infty [$. Thus the function
\begin{equation}
\label{FImplicita}
F\left(\epsilon, \lambda\right)=\frac{J_{\nu_l}\left(\sqrt{\lambda }\right)}{Y_{\nu_l}\left(\sqrt{\lambda }\right)}-\frac{J_{\nu_l}\left(\epsilon \sqrt{\lambda }\right)}{Y_{\nu_l}\left(\epsilon \sqrt{\lambda }\right)}
\end{equation}
is defined for all $(\epsilon ,\lambda )$ in a suitable neighborhood of $(0,\lambda_*)$ in $[0,\infty [\times [0,\infty [$ and is of class $C^1$.
Note that for $\epsilon \ne 0$ the zeros of the function $F$  coincide with the zeros of equation (\ref{EqCaratt}), while the zeros of the function $F(0, \cdot )$
coincide with the zeros of  equation (\ref{EqCarattunp}). We set $\delta  = \epsilon ^{2\nu }  $ and $G(\delta , \lambda )=F(\delta ^{\frac{1}{2\nu }}, \lambda  )$, so that by (\ref{ratio}) it follows that $G$ is a function of class $C^1$ as well. By using standard results (see {\it e.g.,} \cite{olv}), it follows that
\[
\frac{d}
{{dt}}\left( {\frac{{J_{\nu_l}(t)  }}
{{Y_{\nu_l}  \left( t \right)}}} \right) =  - \frac{{2}}
{{\pi tY_{\nu_l} ^2 \left( t\right)}},
\]
hence
\begin{equation}
\label{FLambda}
\frac{\partial G\left(\delta, \lambda\right)}{\partial \lambda} = \frac{1}
{{\pi \lambda }}\left( { - \frac{1}
{{Y_\nu ^2 \left( {\sqrt \lambda  } \right)}} + \frac{1}
{{Y_\nu ^2 \left( {\delta ^{\frac{1}
{{2\nu }}} \sqrt \lambda  } \right)}}} \right)
\end{equation}
and
\begin{equation}
\label{FDeltaX}
\frac{\partial G\left(\delta, \lambda\right)}{\partial \delta}  = \frac{1}
{{\nu \pi \delta Y_\nu ^2 \left( {\delta ^{\frac{1}
{{2\nu }}} \sqrt \lambda  } \right)}} \ \ .
\end{equation}
By passing to the limit in (\ref{FLambda}) and (\ref{FDeltaX}) as $(\delta , \lambda )\to (0,\lambda_*)$ we obtain
\begin{equation}
\label{FLambdaX}
\left .\frac{\partial G\left(\delta, \lambda\right)}{\partial \lambda}\right|_{\left(0, \lambda_*\right)} = - \frac{1}{\pi  \lambda_* } \frac{1}{Y_\nu^2\left(\sqrt{\lambda_* }\right)} \ne 0
\end{equation}
and
\begin{equation}
\label{FDeltaXzero}
\left .\frac{\partial G\left(\delta, \lambda\right)}{\partial \delta}\right|_{\left(0, \lambda_*\right)}   =\frac{\pi }
{{\nu \Gamma ^2 \left( \nu  \right)}}\left( {\frac{\lambda_* }
{4}} \right)^\nu .
\end{equation}
By the Implicit Function Theorem it follows that the zeros of function $G$ in a neighborhood of  $(0, \lambda_*)$ in $[0,\infty [\times [0,\infty [$ are given by the graph of a function $\delta \mapsto \mu (\delta )$ of class $C^1$ such that $\mu (0)=\lambda ^*$ and
\[
\begin{array}{*{20}c}
{\displaystyle \mu\left(\delta\right) = \lambda_* - \left.\frac{G_\delta\left(\delta, \lambda\right)}{G_\lambda\left(\delta, \lambda\right)}\right|_{\left(0,\lambda_*\right)}\delta +o \left(\delta\right),} & {\mathrm{as} \ \delta \to 0^+}.
\end{array}
\]
Taking into account that $|\Omega _{\beta }\setminus \Omega _{\beta }(\epsilon )|=\sigma_{\beta } \epsilon^N/N$ where $\sigma_{\beta }$ denotes the $(N-1)$-dimensional measure of the spherical cap $C_{\beta }$, it follows that the zeros of the function $F$ in a neighborhood of $(0,\lambda_*)$ are given by the graph of a function $\epsilon \mapsto \lambda (\epsilon ) $ such that $\lambda (0)= \lambda _*$ and such that (\ref{main}) is satisfied
with
$$
 a= \frac{\pi ^2 N^{\frac{2\nu_l}{N}} \lambda_*^{\nu_l+1}Y_{\nu_l}^2(\sqrt{\lambda_*})  }{\nu_l4^{\nu_l}  \sigma_{\beta }^{\frac{2\nu_l}{N}}  \Gamma^2(\nu_l ) }\, .
$$
Formula (\ref{main2}) follows by (\ref{main}) and by recalling that if $\lambda_*=\lambda_1[\Omega_{\beta} ]$  then $\lambda_*$ is a solution to  equation (\ref{EqCarattunp}) with  $l=l_{\beta}$ (see the proof of Theorem~\ref{thmsum}). \hfill $\Box$\\

\begin{rem}\label{mainrem} By Theorem~\ref{thmsum}, if $\beta \in ]\pi /2 , \pi [$ then the eigenfunctions of the Dirichlet Laplacian on $\Omega _{\beta }$ satisfy condition (\ref{sumcond}) for all $p\in [2, N/(1-l_{\beta })[$. The limiting case $p= N/(1-l_{\beta })$ has to be excluded. However, we note that  if $p= N/(1-l_{\beta })$ then
$$
1-\frac{2}{p}=  \frac{N+2l_{\beta } -2}{N}
$$
which is exactly the exponent of the leading term  in (\ref{main2}). As it is proved in Theorem~\ref{shathm}, this implies
that in estimate (\ref{bulaest}) one cannot replace the exponent $1-2/p$ by a better exponent of the type $1-2/p +\delta (p)$ where $\delta (p)>0$ depends with continuity on $p$.
\end{rem}

\section{Sharpness of the stability estimate in open sets of class $C^{0,1}_M({\mathcal{A}})$}

The stability estimates in \cite{buladir} are proved for open sets of class $C^{0,1}_M({\mathcal{A}})$. We recall the definition of this class.

For any set $V$ in ${\mathbb{R}}^N$ and $\delta >0$ we denote by $V_{\delta }$ the set $\{x\in V:\ d(x, \partial V )>\delta \}$. Moreover, by a rotation in ${\mathbb{R}}^N$ we mean a $N\times N$-orthogonal matrix with real entries which we identify with the corresponding linear operator acting in ${\mathbb{R}}^N$.\\

\begin{defn}
\label{class}

Let $ \rho >0$, $s,s'\in\mathbb{N}$, $s'< s$
and  $\{V_j\}_{j=1}^s$ be a family of bounded open cuboids  and
$\{r_j\}_{j=1}^{s} $ be a family of rotations in ${\mathbb{R}}^N $.

We say that that ${\mathcal{A}}= (  \rho , s,s', \{V_j\}_{j=1}^s, \{r_j\}_{j=1}^{s} ) $ is an atlas in ${\mathbb{R}}^N$ with the parameters
$\rho , s,s', \{V_j\}_{j=1}^s, \{r_j\}_{j=1}^{s}$, briefly an atlas in ${\mathbb{R}}^N$.

Let $M>0$. We denote by $C^{0,1}_M( {\mathcal{A}}   )$ the family of all open sets $\Omega $ in ${\mathbb{R}}^N$
satisfying the following properties:

(i) $ \Omega\subset \bigcup\limits_{j=1}^s(V_j)_{\rho}$ and $(V_j)_\rho\cap\Omega\ne\emptyset;$

(ii) $V_j\cap\partial \Omega\ne\emptyset$ for $j=1,\dots s'$, $ V_j\cap \partial\Omega =\emptyset$ for $s'<j\le s$;

(iii) for $j=1,...,s$
$$
r_j(V_j)=\{\,x\in \mathbb{R}^N:~a_{ij}<x_i<b_{ij}, \,i=1,....,N\},
$$

\noindent and

$$
r_j(\Omega\cap V_j)=\{x\in\mathbb{R}^N:~a_{1j}<x_1<g_j(\bar x),~\bar x\in W_j\},$$

\noindent where $\bar x=(x_2,...,x_{N})$, $W_j=\{\bar
x\in\mathbb{R}^{N-1}:~a_{ij}<x_i<b_{ij},\,i=2,...,N\}$
and $g_j$ is a continuous real-valued function defined on $\overline {W}_j$ (it is meant that if $s'<j\le s$ then $g_j(\bar x)=b_{1j}$ for all $\bar x\in \overline{W}_j$);

moreover for $j=1,\dots ,s'$
$$
a_{1j}+\rho\le g_j(\bar x)\le b_{1j}-\rho ,$$

\noindent for all $\bar x\in \overline{W}_j$ and the Lipschitz constant of $g_j$ satisfies
$$
{\rm Lip }\, g_j \le M.
$$

We say that an open set is of class $C^{0,1}$ if there exists an atlas ${\mathcal{A}}$ and $M>0$ such that $\Omega$ is of class $C^{0,1}_M({\mathcal{A}})$.
\end{defn}

In Theorem~\ref{shathm} we prove that the following result is sharp (see \cite[Cor.~5.16]{buladir} for the original statement concerning a general class
of second order uniformly elliptic operators).

\begin{thm}\label{bulathm}
Let ${\mathcal{A}}$ be an atlas in ${\mathbb{R}}^N$ and $M>0$. Let $\Omega _1$ be an open set in ${\mathbb{R}}^N$ of class $C^{0,1}_M({\mathcal{A}})$.
Assume that there exists $p\in ]2,\infty ] $ such that condition (\ref{sumcond})
is satisfied for  all eigenfunctions $u $ of the Dirichlet Laplacian on $\Omega_1$.
Then for each $n\in {\mathbb{N}}$ there exists $c_n>0 $ such that estimate (\ref{bulaest}) holds
for all $\Omega_2\in C^{0,1}_M({\mathcal{A}})$ such that $\Omega_2\subset \Omega_1$ and $|\Omega_1\setminus \Omega_2|<c_n^{-1}$.
\end{thm}

Note that by combining Theorem~\ref{bulathm} and Jerison and Kenig~\cite[Theorems~0.5, 1.1, 1.3]{jeke} one can immediately deduce the following statement.

\begin{corol}\label{dimecor}  Let ${\mathcal{A}}$ be an atlas in ${\mathbb{R}}^N$ and $M>0$. Let $\Omega _1$ be an open set in ${\mathbb{R}}^N$ of class $C^{0,1}_M({\mathcal{A}})$. Then there exists $r>0$ depending only on $M$ with $r>1/2$ if $N=2$ and $r>1/3$ if $N\geq 3$ such that
\begin{equation}
\label{bulaesterre}
|\lambda _n[\Omega _1]- \lambda _n[\Omega _2]|\le c_n |\Omega _1\setminus \Omega _2|^r,
\end{equation}
for all $\Omega_2\in C^{0,1}_M({\mathcal{A}})$ such that $\Omega_2\subset \Omega_1$ and $|\Omega_1\setminus \Omega_2|<c_n^{-1}$.
\end{corol}

In Theorem~\ref{shathm} we shall also prove that for $N=2,3$, in Corollary~\ref{dimecor} one cannot choose $r$ to be independent of $M$.

It is clear that for each $\epsilon \in [0,1[$ the sets
$ \Omega_{\beta }(\epsilon )$ defined in (\ref{omegaeps}) are of class $C^{0,1}$. However, if $\beta >\pi /2$ there is not a fixed  atlas ${\mathcal{A}}$  such that
$\Omega _{\beta }(\epsilon )$ is of class $C^{0,1 }_M({\mathcal{A}})$ for all $\epsilon \in [0,1[$.
Thus, in order to prove that the result of Theorem~\ref{bulathm} is sharp, we cannot use directly the open sets $\Omega_{\beta }(\epsilon )$.
To overcome this difficulty, we need to modify the sets $\Omega_{\beta }(\epsilon )$ and replace them by suitable open sets $\tilde \Omega _{\beta }(\epsilon )$  belonging to a class defined by the same fixed atlas.

For all $\epsilon \in [0,1 [$ we set
$$
\tilde \Omega _{\beta }(\epsilon )= \left\{x=(x_1,\bar x )\in {\mathbb{R}}^N:\ g(\bar x )<x_1,\ \ |x |<1  \right\}
$$
where $\bar x=(x_2, \dots ,x_N)$,  $g(\bar x)= \epsilon - |\bar x|\tan \frac{\beta }{2}$ if $|\bar x|\le \epsilon \sin \beta $, and $g(\bar x)=|\bar x| \cot \beta $ if
$|\bar x|> \epsilon \sin \beta $.

We note that
\begin{equation}\label{inclus}
\Omega _{\beta }(\epsilon )\subset  \tilde \Omega_{\beta } (\epsilon )\subset \Omega_{\beta }(A\epsilon )\subset \Omega_{\beta } ,
\end{equation}
for all $\epsilon \in [0,1[$,  where $A=\sin \beta / \sqrt { 2(1-\cos \beta ) }$ (see the figure).

The first inclusion in (\ref{inclus}) can be proved by observing that  if  $(x_1, \bar x )\in \Omega_{\beta }(\epsilon )$ then $x_1= r\cos \theta _1$ and $|\bar x|= r \sin \theta _1$ where $\epsilon <r<1$  and $\theta _1$ is as in (\ref{omegaN2}), (\ref{omegaN34}). Then one can easily verify that $x_1>g(\bar x )$ for all $(x_1,\bar x )\in \Omega _{\beta }(\epsilon )$ which implies that $(x_1, \bar x)\in  \tilde\Omega_{\beta }(\epsilon )$.  The second inclusion in (\ref{inclus}) can be proved in a similar way. Indeed, by means of a simple computation one can prove that  if $x=(x_1,  \bar x )\in  \tilde \Omega_{\beta }(\epsilon )$ then $r=|x|> \epsilon \sin \beta /\sqrt { 2(1-\cos \beta ) }  $; moreover, the spherical coordinate  $\theta_1 $ of the point $x$ satisfies the appropriate inequalities in (\ref{omegaN2}), (\ref{omegaN34}), since $\cot \theta _1 = x_1/ |\bar x|  $; thus $x\in  \Omega_{\beta }(A\epsilon )$.\\
Note that in particular $\tilde \Omega_{\beta }(0)=\Omega_{\beta }$.

\begin{figure}[hb]
\centering
\includegraphics[angle=270,width=8cm,viewport=0 0 600 650,clip]{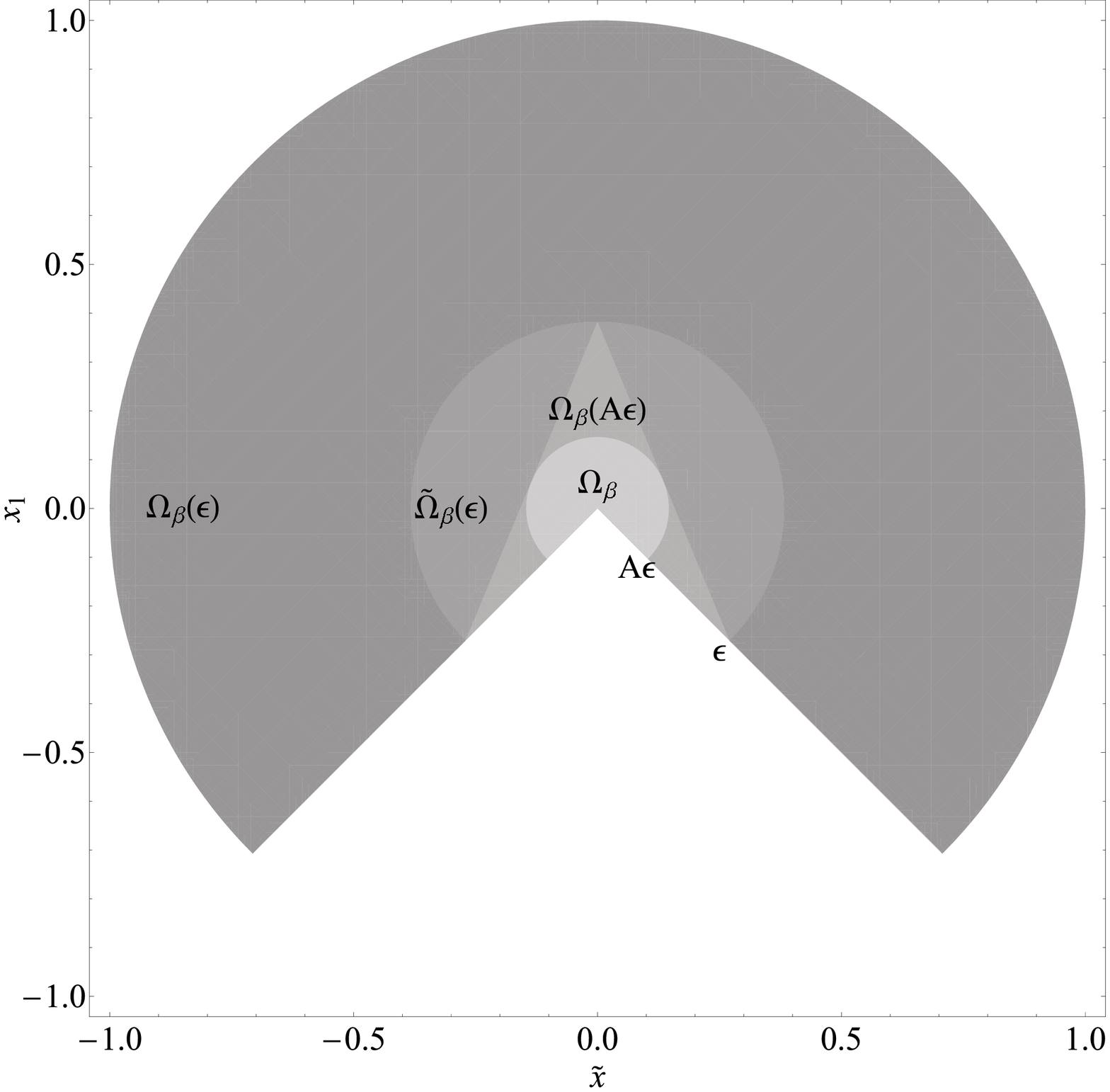}
\end{figure}

Also, it is clear that there exists an atlas ${\mathcal{A}}$ and $M>0$ such that
$$
\tilde \Omega _{\beta }(\epsilon )\ \ {\rm is}\ {\rm of}\ {\rm class }\ \ C^{0,1}_M({\mathcal {A}}),
$$
for all $\epsilon \in [0, 1/2 [$.

\begin{lem}
Let $\beta \in ]0,\pi [$. There exist $d_1,d_2>0$ such that
\begin{equation}
\label{esttilde}
d_1|\Omega_{\beta }\setminus \tilde \Omega_{\beta }(\epsilon ) | ^{\frac{N+2l_{\beta } -2}{N}  } \le
|\lambda_1[\tilde \Omega_{\beta }(\epsilon )]-\lambda_1[ \Omega_{\beta }]|\le d_2|\Omega_{\beta }\setminus \tilde \Omega_{\beta }(\epsilon ) |  ^{\frac{N+2l_{\beta } -2}{N}},
\end{equation}
for all $\epsilon >0$ sufficiently small.
\end{lem}

{\bf Proof.} By the monotonicity of the eigenvalues of the Dirichlet Laplacian and inclusions (\ref{inclus}), it follows that
\begin{equation}\lambda_1[\Omega_{\beta }]\le
\lambda_1[\Omega_{\beta }(A\epsilon )]\le \lambda_1[\tilde\Omega_{\beta }(\epsilon )] \le \lambda_{1}[\Omega _{\beta }(\epsilon )].
\end{equation}
Thus, by Theorem~\ref{mainthm} it follows that there exist $\tilde d_1, \tilde d_2>0$ such that
$$
\tilde d_1|\Omega_{\beta }\setminus \Omega_{\beta }(A\epsilon )  |^{\frac{N+2l_{\beta } -2}{N}}\le |\lambda_1[\tilde \Omega_{\beta }(\epsilon )]-\lambda_1[ \Omega_{\beta }]|  \le \tilde d_2|  \Omega _{\beta }\setminus \Omega_{\beta }(\epsilon )|^{\frac{N+2l_{\beta } -2}{N}} .
$$
Inequality (\ref{esttilde}) follows by noting that $|\Omega _{\beta }\setminus \tilde \Omega_{\beta }(A\epsilon )|\le  |   \Omega _{\beta }\setminus \Omega_{\beta }(A\epsilon )|$, $|\Omega_{\beta }\setminus \Omega_{\beta }(\epsilon )|\le |  \Omega _{\beta }\setminus \tilde \Omega_{\beta }(\epsilon /A)|$ and $|\Omega _{\beta }\setminus \tilde \Omega_{\beta }(c\epsilon )|=c^N|\Omega _{\beta }\setminus \tilde \Omega_{\beta }(\epsilon )|$ for all $c\in ]0,1]$.\hfill $\Box $\\

Finally, we can prove the main result of this section.

\begin{thm}\label{shathm}
For each $N\geq 2$ the exponent $1-2/p $ in the right-hand side of estimate (\ref{bulaest}) in Theorem~\ref{bulathm} cannot be replaced by $1-2/p+\delta (p) $ where  $\delta (p) >0$ depends with continuity on  $p$.

Moreover, if  $N=2,3$  then for any fixed $r>1/N$ there exists an open set $\Omega _1$ of class $C^{0,1}$ such that estimate (\ref{bulaesterre}) does not hold.
\end{thm}

{\bf Proof.} Let $\beta \in ]\pi /2,\pi  [$. By Lemma~\ref{charcap},  $l_{\beta }<1$. Let  $\Omega_1=\Omega_{\beta }$ and $\Omega_2=\tilde\Omega_{\beta }(\epsilon )$ for $\epsilon \in ]0,1/2[$.
Recall that $\Omega_2\subset \Omega_1$ and that $\Omega_1$, $\Omega_2$  are open sets of the same class $C^{0,1}_M({\mathcal{A}})$ for a suitable atlas ${\mathcal{A}}$ and $M >0$ independent of $\epsilon$.

By Theorem~\ref{thmsum} the eigenfunctions of the Dirichlet Laplacian on $\Omega_1$ satisfy condition (\ref{sumcond}) for all  $p\in [2, N/(1-l_{\beta })[$. Moreover,
\begin{equation}
\label{limite}
\lim_{p\to  N/(1-l_{\beta })}1-\frac{2}{p}= \frac{N+2l_{\beta } -2}{N}.
\end{equation}

Assume now by contradiction that the exponent $1-2/p $ in the right-hand side of (\ref{bulaest}) in Theorem~\ref{bulathm} could be replaced by $1-2/p+\delta (p) $ where $\delta $ is as in the statement. Then by (\ref{limite}) one could find $\bar p$ sufficiently close to $N/(1-l_{\beta }) $ such that
$$
1-\frac{2}{\bar p}+\delta (\bar p) > \frac{N+2l_{\beta } -2}{N}.
$$
Thus by applying Theorem \ref{bulathm} with the exponent $1-2/\bar p +\delta (\bar p )  $  in (\ref{bulaest}), we would obtain that the second inequality  in (\ref{esttilde}) holds with an exponent larger than $\frac{N+2l_{\beta } -2}{N}$, which contradicts the validity of the first inequality  in (\ref{esttilde}) for $\epsilon $
sufficiently small.

The second part of the statement follows in a similar way by choosing $\Omega_1=\Omega_{\beta }$ and by observing that by Remark~\ref{camerarem}
$$\lim_{\beta\to \pi}\frac{N+2l_{\beta } -2}{N}=\frac{1}{N}$$
whenever $N=2,3$.
\hfill $\Box $\\

{\bf Acknowledgments.} The authors are thankful to Professor Victor I. Burenkov for encouraging this investigation and useful discussions.
The authors would also like to thank Professor Giuseppe Savar\'{e} for pointing out reference \cite{jeke}.

\vspace{1cm}
\noindent
{\small Pier Domenico Lamberti and Marco Perin\\
Dipartimento di Matematica Pura ed Applicata\\
Universit\`{a} degli Studi di Padova\\
Via Trieste 63\\
35121 Padova\\
Italy\\
e-mail: lamberti@math.unipd.it, marco\_perin1765@alice.it
}
\end{document}